\documentclass[11pt, a4paper,reqno]{amsart}
\usepackage{amsmath}
\usepackage{amsfonts}
\usepackage{amssymb}
\usepackage{amsthm}

\numberwithin{equation}{section}

\newtheorem{theorem}[equation]{Theorem}
\newtheorem{lemma}[equation]{Lemma}

\theoremstyle{definition}
\newtheorem{definition}[equation]{Definition}

\theoremstyle{remark}
\newtheorem{remark}[equation]{Remark}

\def\kint_#1{\mathchoice%
          {\mathop{\kern 0.2em\vrule width 0.6em height 0.69678ex depth -0.58065ex
                  \kern -0.8em \intop}\nolimits_{\kern -0.4em#1}}%
          {\mathop{\kern 0.1em\vrule width 0.5em height 0.69678ex depth -0.60387ex
                  \kern -0.6em \intop}\nolimits_{#1}}%
          {\mathop{\kern 0.1em\vrule width 0.5em height 0.69678ex depth -0.60387ex
                  \kern -0.6em \intop}\nolimits_{#1}}%
          {\mathop{\kern 0.1em\vrule width 0.5em height 0.69678ex depth -0.60387ex
                  \kern -0.6em \intop}\nolimits_{#1}}}
\def\vintslides_#1{\mathchoice%
          {\mathop{\kern 0.1em\vrule width 0.5em height 0.697ex depth -0.581ex
                  \kern -0.6em \intop}\nolimits_{\kern -0.4em#1}}%
          {\mathop{\kern 0.1em\vrule width 0.3em height 0.697ex depth -0.604ex
                  \kern -0.4em \intop}\nolimits_{#1}}%
          {\mathop{\kern 0.1em\vrule width 0.3em height 0.697ex depth -0.604ex
                  \kern -0.4em \intop}\nolimits_{#1}}%
          {\mathop{\kern 0.1em\vrule width 0.3em height 0.697ex depth -0.604ex
                  \kern -0.4em \intop}\nolimits_{#1}}}

\newcommand{\eps}{\varepsilon}

\newcommand{\Rn}{\mathbb{R}^d}

\newcommand{\esssup}{\operatornamewithlimits{ess\, sup}}
\newcommand{\essinf}{\operatornamewithlimits{ess\,inf}}
\newcommand{\essosc}{\operatornamewithlimits{ess\, osc}}

\renewcommand{\div}{\nabla \cdot}

\renewcommand{\l}{\left}
\renewcommand{\r}{\right}

\def\Xint#1{\mathchoice
{\XXint\displaystyle\textstyle{#1}}%
{\XXint\textstyle\scriptstyle{#1}}%
{\XXint\scriptstyle\scriptscriptstyle{#1}}%
{\XXint\scriptscriptstyle\scriptscriptstyle{#1}}%
\!\int}

\def\XXint#1#2#3{{\setbox0=\hbox{$#1{#2#3}{\int}$}
\vcenter{\hbox{$#2#3$}}\kern-.5\wd0}}

\def\dashint{\Xint-}

\title[Local H\"older continuity]{Local H\"older continuity for doubly nonlinear parabolic equations}
\author[Kuusi, Siljander and Urbano]{Tuomo Kuusi, Juhana Siljander and Jos\'e Miguel Urbano}

\begin{document}

\begin{abstract}
We give a proof of the H\"older continuity of weak solutions of certain degenerate doubly nonlinear parabolic equations in measure spaces. We only assume the measure to be a doubling non-trivial Borel measure which supports a Poincar\'e inequality. The proof discriminates between large scales, for which a Harnack inequality is used, and small scales, that require intrinsic scaling methods.
\end{abstract}



\thanks{Research of JMU supported by CMUC/FCT and project UTAustin/MAT/0035/2008.}

\maketitle

\section{Introduction}
We consider the regularity issue for nonnegative weak solutions of the
doubly nonlinear parabolic equation
\begin{equation}\label{equation}
\frac{\partial (u^{p-1})}{\partial t}-\div{(|\nabla u|^{p-2}\nabla u)}=0, \quad 2\le p<\infty.
\end{equation}
This equation is a prototype of a parabolic equation of $p$-Laplacian type. Its solutions can be scaled by nonnegative
factors, but in general we cannot add a constant to a solution so that the resulting function
would be a solution to the same equation. 

The purpose of this paper is to obtain a clear and transparent proof for the local H\"older continuity of nonnegative weak solutions of \eqref{equation}. Our work is a continuation to~\cite{KinnKuus07}, where Harnack's inequality for the same equation is proved. See also \cite{Trud68}, \cite{GianVesp06c}, \cite{GianVesp06a} and \cite{Vesp94}. However, since we cannot add constants to solutions, the Harnack estimates do not directly imply the local H\"older continuity. To show that our proof is based on a general principle, we consider the case where the Lebesgue measure is replaced by a more general Borel measure, which is merely assumed to be doubling and to support
a Poincar\'e inequality. In the weighted case, parabolic equations have earlier been studied in
\cite{ChiaSera84a}, \cite{ChiaSera84b} and \cite{Surn10}. See also \cite{FabeKeniSera82}.

This kind of doubly nonlinear equations have been considered by Vespri \cite{Vesp92}, Porzio and Vespri \cite{PorzVesp93}, and Ivanov \cite{Ivan91}, \cite{Ivan95}. The known regularity proofs are based on the method of intrinsic scaling, originally introduced by DiBenedetto, and they seem to depend highly on the
particular form of the equation. However, the passage from one equation to another is not completely clear. For other parabolic equations, the problem has been studied at length, see \cite{DiBe93}, \cite{DiBe86}, \cite{DiBeUrbaVesp04} and \cite{Urba08}, and the references therein.



The difficulty with equation \eqref{equation} is that there is a certain kind of dichotomy in its behavior. Correspondingly, the proof has been divided in two complementary cases:
\[
 \text{Case I}:\quad 0\le\essinf{u} <<\essosc{u}
\]
and
\[
 \text{Case II}:\quad u^{p-2}u_t \approx Cu_t. \ \ \ \ \ \ \ \  \ \ \ \ \ \ \
\]

In large scales, i.e., in Case I, the scaling property of the equation dominates and, consequently, the reduction of the oscillation follows immediately from Harnack's inequality. In small scales, on the other hand, the oscillation is already very small and thus the solution itself is between two constants, the infimum and the supremum, whose difference is negligible. Correspondingly, the nonlinear time derivative term, which formally looks like $u^{p-2}u_t$, behaves like $u_t$ and we end up with a $p$-parabolic type behavior. However, also in this case, we still need to modify the known arguments. In particular, the energy estimates are not available in the usual form and we need to use modified versions as in \cite{DiBeFrie85a}, \cite{Ivan95} and \cite{Zhou94}.

Our argument also applies to doubly nonlinear equations of $p$-Laplacian type that are of the form
\[
\frac{\partial (u^{p-1})}{\partial t}- \div\mathcal{A}(x,t,u,\nabla u) = 0,
\]
with $\mathcal{A}(x,t,\cdot,\cdot)$ satisfying the usual structure assumptions.
For expository purposes, we only consider~\eqref{equation}.


Very recently, a direct geometric method to obtain local H\"older continuity for parabolic equations has
been developed in \cite{DiBeGianVesp10} and \cite{GianSurnVesp10}.
Despite the effort, the general picture remains unclear.

\section{Preliminaries} \label{section:preliminaries}

Let $\mu$ be a Borel measure and $\Omega$ be an open set in $\Rn$.
The Sobolev space $H^{1,p}(\Omega;\mu)$ is defined to be the completion of $C^\infty(\Omega)$ with respect to the Sobolev norm
\[
\|u\|_{1,p,\Omega}=\left(\int_\Omega(|u|^p + |\nabla u|^p)\, d\mu\right)^{1/p}.
\]
A function $u$ belongs to the local Sobolev space $H_{loc}^{1,p}(\Omega;\mu)$ if it belongs to
$H^{1,p}(\Omega';\mu)$ for every $\Omega' \Subset \Omega$. Moreover, the Sobolev space with zero boundary values $H^{1,p}_0(\Omega;\mu)$ is defined as the completion of $C_0^\infty(\Omega)$ with respect to the Sobolev norm. For more properties of Sobolev spaces, see e.g. \cite{HeinKilpMart93}.

Let $t_1<t_2$.
The parabolic Sobolev space $L^p(t_1,t_2;H^{1,p}(\Omega;\mu))$ is the space of functions $u(x,t)$
such that, for almost every $t$, with $t_1<t<t_2$, the function $u(\cdot, t)$ belongs to $H^{1,p}(\Omega;\mu)$ and
\[
\int_{t_1}^{t_2}\int_\Omega (|u|^p+|\nabla u|^p) \, d\nu < \infty,
\]
where we denote $d\nu=d\mu\,dt$.

The definition of the space $L_{loc}^p(t_1,t_2;H_{loc}^{1,p}(\Omega;\mu))$ is analogous.

\begin{definition}
A function $u \in L_{loc}^p(t_1,t_2;H_{loc}^{1,p}(\Omega;\mu))$ is a
weak solution of equation \eqref{equation} in $\Omega\times(t_1,t_2)$ if it satisfies the integral equality
\begin{equation} \label{weak_solution}
\int_{t_1}^{t_2}\int_{\Omega} \left( |\nabla u|^{p-2}\nabla u\cdot \nabla \phi
-u^{p-1} \frac{\partial\phi}{\partial t} \right)\, d\nu = 0
\end{equation}
for every  $\phi \in C_0^\infty(\Omega \times (t_1,t_2))$.
\end{definition}


Next, we recall a few definitions and results from analysis on metric measure spaces.
The measure $\mu$ is doubling if there is a universal constant $D_0\ge 1$ such that
\[
 \mu(B(x,2r))\le D_0 \mu(B(x,r)),
\]
for every $B(x,2r)\subset\Omega$. Here
\[
B(x,r)=\{y\in \Rn : |y-x|<r\}
\]
denotes the standard open ball in $\Rn$.
Let $0<r<R<\infty$.
A simple iteration of the doubling condition implies that
\[
\frac{\mu(B(x,R))}{\mu(B(x,r))}
\le C\left(\frac Rr\right)^{d_\mu},
\]
where $d_\mu=\log_2 D_0$.
A doubling measure in $\Omega$ also satisfies the following annular decay property.
There exist constants $0<\alpha<1$ and $c\ge 1$ such that
\begin{equation}\label{annular_decay}
\mu(B(x,r)\setminus B(x,(1-\delta)r))\le c\delta^\alpha\mu(B(x,r)),
\end{equation}
for all $B(x,r)\subset\Omega$ and $0<\delta<1$.

The measure is said to support a weak $(1,p)$-Poincar\'e inequality if there exist constants $P_0>0$ and $\tau\ge 1$ such that
\[
\dashint_{B(x,r)}|u-u_{B(x,r)}| \, d\mu
\le P_0 r\left(\dashint_{B(x,\tau r)} |\nabla u|^p \, d\mu\right)^{1/p},
\]
for every $u\in H^{1,p}_{loc}(\Omega;\mu)$ and $B(x,\tau r)\subset\Omega$.
Here, we denote
\[
 u_{B(x,r)}=\dashint_{B(x,r)} u \, d\mu = \frac{1}{\mu(B(x,r))}\int_{B(x,r)} u\, d\mu.
\]
The word weak refers to the constant $\tau$, that may be strictly greater than one.
In $\Rn$ with a doubling measure, the weak $(1,p)$-Poincar\'e inequality
with some $\tau\ge1$ implies the $(1,p)$-Poincar\'e inequality with $\tau=1$,
see Theorem 3.4 in \cite{HajlKosk00}.
Hence, we may assume that $\tau=1$.

On the other hand, the weak $(1,p)$-Poincar\'e inequality and the doubling condition
imply a weak $(\kappa,p)$-Sobolev-Poincar\'e inequality with
\begin{equation} \label{kappa}
\kappa =
\begin{cases}
\dfrac{d_\mu p}{d_\mu -p}, & 1<p< d_\mu, \\
2p, & p \ge d_\mu,
\end{cases}
\end{equation}
where $d_\mu$ is as above.
In other words, Poincar\'e and doubling imply the Sobolev inequalities.
More precisely, there are constants $C>0$ and $\tau'\ge1$ such that
\begin{equation}\label{Sobolev}
\Big(\dashint_{B(x,r)}|u-u_{B(x,r)}|^\kappa\,d\mu\Big)^{1/\kappa}
\le Cr\left(\dashint_{B(x,\tau' r)}|\nabla u|^p\,d\mu\right)^{1/p},
\end{equation}
for every $B(x,\tau'r)\in\Omega$.
The constant $C$ depends only on $p$, $D_0$ and $P_0$.
For the proof, we refer to \cite{HajlKosk00}.
Again, by Theorem 3.4 in \cite{HajlKosk00} we may take $\tau'=1$ in \eqref{Sobolev}.

For Sobolev functions with the zero boundary values,
we have the following version of Sobolev's inequality.
Suppose that $u\in H_0^{1,p}(B(x,r);\mu)$.
Then
\begin{equation}\label{Sobolevzero}
\left(\dashint_{B(x,r)} |u|^\kappa\, d \mu
\right)^{1/\kappa} \le C r \left(\dashint_{B(x,r)}
|\nabla u|^p\, d \mu \right)^{1/p}.
\end{equation}
For the proof we refer, for example, to \cite{KinnShan01}.

Moreover, by a recent result in \cite{KeitZhon08},
the weak $(1,p)$-Poincar\'e inequality and the doubling condition also imply
the $(1,q)$-Poincar\'e inequality for some $q<p$, that is
\begin{equation}\label{poincare}
\dashint_{B(x,r)}|u-u_{B(x,r)}|\, d\mu
\le C r\left(\dashint_{B(x,r)} |\nabla u|^q \, d\mu\right)^{1/q}.
\end{equation}
Consequently, also \eqref{Sobolev} holds with $p$ replaced by $q$. We also obtain the $(q,q)$-Poincar\'e inequality for some $q<p$.



In the sequel, we shall refer to \textit{data} as the set of a priori constants
$p$, $d$, $D_0$, and $P_0$.

Our main result is the following theorem.

\begin{theorem}\label{main_theorem}
Let $2\le p<\infty$ and assume that the measure is doubling, supports a weak $(1,p)$-Poincar\'e inequality and is non-trivial in the sense that the measure of every non-empty open set is strictly positive and the measure of every bounded set is finite. Moreover, let $u\ge 0$ be a weak solution of equation~\eqref{equation} in $\Rn$. Then $u$ is locally H\"older continuous.
\end{theorem}

We will use the following notation for balls and cylinders, respectively:
\[
B(r)=B(0,r)
\]
and
\[
Q_t(s,r)=B(r)\times(t-s,t).
\]
For simplicity, we will also denote
\[
Q(s,r)=Q_0(s,r)=B(r)\times (-s,0).
\]

Recall Harnack's inequality from \cite{KinnKuus07}.

\begin{theorem}\label{Harnack}
Let $1<p<\infty$ and suppose that the measure $\mu$ is doubling and supports a weak $(1,p)$-Poincar\'e inequality.
Moreover, let $u\ge 0$ be a weak solution to ~\eqref{equation} in $\Rn$.
Then there exists a constant $H_0=H_0(p,d, D_0, P_0, (t-(s-r^p))/r^p)\ge 2$ such that
\[
 \esssup_{Q_t(r^p,r)}{u}\le H_0 \essinf_{Q_{s}(r^p,r)}{u},
\]
where $s> t+r^p$.
\begin{proof}
See~\cite{KinnKuus07}.
\end{proof}
\end{theorem}
In addition, in \cite{KinnKuus07} it is also proved that all solutions of equation~\eqref{equation} are locally bounded. In the sequel, we will assume this knowledge without any further comments.

\section{Constructing the setting}

Our proof is based on the known classical argument of reducing the oscillation, see \cite{DiBe93}, \cite{DiBeUrbaVesp04} and \cite{Urba08}. However, the equation under study has some intrinsic properties which are not present, for instance, in the case of the $p$-parabolic equation. In large scales, the scaling property dominates and the oscillation reduction follows easily from Harnack's inequality. In this case, the equation resembles the usual heat equation.

In small scales, in turn, the equation changes its behavior to look more like the evolution $p$-Laplace equation. Indeed, when we zoom in by reducing the oscillation, the infimum and the supremum get closer and closer to each other. Consequently, the weight $u^{p-2}$ in the time derivative term starts to behave like a constant coefficient and we end up with a $p$-parabolic type behavior. Resembling this divide between large and small scales, the proof has to be divided in two cases. 

We study the (local) H\"older continuity in a compact set $K$ and we choose the following numbers accordingly. Let
\[
\mu_0^-\le\essinf_Ku\quad\text{and}\quad\mu_0^+\ge\esssup_Ku,
\]
and define
\[
 \omega_0=\mu_0^+-\mu_0^-.
\]
Furthermore, choose $\mu_0^-$ small enough so that
\begin{equation}\label{bounds_inf}
(2H_0+1)\mu_0^-\le \omega_0
\end{equation}
holds. We will construct an increasing sequence $\{\mu_i^-\}$ and a decreasing sequence $\{\mu_i^+\}$ such that
\[
 \mu_i^+-\mu_i^-=\omega_i=\sigma^i\omega_0
\]
for some $0<\sigma<1$. Moreover, these sequences can be chosen so that
\[
 \esssup_{Q^i}u\le \mu_i^+
\]
and
\[
\essinf_{Q^i}u\ge \mu_i^- ,
\]
for some suitable sequence $\{Q^i\}$ of cubes. Consequently,
\[
\essosc_{Q^i}u\le \omega_i.
\]
The cubes here will be chosen so that their size decreases in a controllable way, from which we can deduce the H\"older continuity. Observe also that if
\[
(2H_0+1)\mu_{j_0}^- \leq \omega_{j_0}
\]
 fails for some $j_0$, the above sequences have been chosen so that
\begin{equation}\label{elliptic_harnack}
 \frac{\mu_j^+}{\mu_j^-}< 2H_0+2
\end{equation}
for all $j\ge j_0$.

%

We are studying the local H\"older continuity in a compact set $K$. Our aim is to show that the oscillation around any point in $K$ reduces whenever we suitably decrease the size of the set where the oscillation is studied. 

The next step is to iterate this reduction process. We end up with a sequence of cylinders $Q^i$. For all purposes, in the sequel, it is enough to study the cylinder $Q^0:=Q(\eta r^p,r)$ instead of the set $K$. Indeed, for any point in $K$ we can build the sequence of suitable cylinders, but since we can always translate the equation, we can, without loss of generality, restrict the study to the origin.

The equation~\eqref{equation} has its own time geometry too, that we need to respect in the arguments. This is important, especially in small scales, when the equation resembles the evolution $p$-Laplace equation. We will use a scaling factor $\eta=2^{\lambda_1(p-2)+1}$
in the time direction, where $\lambda_1\ge 1$ is an a priori constant to be determined later.

\section{Fundamental estimates}

We start the proof of Theorem~\ref{main_theorem} by proving the usual energy estimate in a slightly modified setting,
which overcomes the problem that we cannot add constants to solutions, see \cite{DiBeFrie85a}, \cite{Ivan95} and \cite{Zhou94}.
We introduce the auxiliary function
\begin{align*}
\mathcal{J}((u{-}k)_\pm)= & \pm \int_{k^{p-1}}^{u^{p-1}} \left(\xi^{1/(p-1)} {-} k \right)_\pm \, d\xi 
\\ = & \pm (p-1)\int_{k}^{u} \left(\xi{-}k \right)_\pm \xi^{p-2} \, d\xi 
\\ = & (p-1)\int_0^{(u{-}k)_\pm}(k\pm \xi)^{p-2}\xi \, d\xi.
\end{align*}
Hence, we have
\begin{equation}\label{J_derivative}
\frac{\partial }{\partial t}\mathcal{J}((u{-}k)_\pm)
=\pm\frac{\partial (u^{p-1})}{\partial t}(u{-}k)_\pm.
\end{equation}

In the sequel, we will also need the following estimates. Clearly,
\begin{equation}\label{plus_upper_estimate}
\begin{split}
\mathcal{J}((u{-}k)_+)&=(p-1)\int_0^{(u{-}k)_+}(k+ \xi)^{p-2}\xi \, d\xi \\
&\le\frac{p-1}{2}(k+(u{-}k)_+)^{p-2}(u{-}k)_+^2 \\
&\le\frac{p-1}{2}u^{p-2}(u{-}k)_+^2 \\
\end{split}
\end{equation}
and
\begin{equation}\label{plus_lower_estimate}
\begin{split}
\mathcal{J}((u{-}k)_+)
&\ge (p-1)k^{p-2}\int_0^{(u{-}k)_+}\xi\, d\xi \\
&\ge (p-1)k^{p-2}\frac{(u{-}k)_+^2}{2}.
\end{split}
\end{equation}
Observe that the assumption $p\ge 2$ is used here.

On the other hand,
\begin{equation}\label{minus_lower_estimate2}
\begin{split}
\mathcal{J}((u{-}k)_-)&=(p-1)\int_0^{(u{-}k)_-}(k- \xi)^{p-2}\xi \, d\xi \\
&\ge \frac{(p-1)}{2}u^{p-2}(u{-}k)_-^2.
\end{split}
\end{equation}
Moreover,
\begin{equation}\label{minus_upper_estimate}
\begin{split}
\mathcal{J}((u{-}k)_-)&=(p-1)\int_0^{(u{-}k)_-}(k- \xi)^{p-2}\xi \, d\xi \\
&\le(p-1)k^{p-2}\int_0^{(u{-}k)_-}\xi\, d\xi \\
&=(p-1)k^{p-2}\frac{(u{-}k)_-^2}{2}.
\end{split}
\end{equation}

Now we are ready for the fundamental energy estimate.

\begin{lemma}\label{energy}
Let $u\ge 0$ be a weak solution of \eqref{equation} and let $k\ge 0$.
Then there exists a constant $C=C(p)>0$ such that
\begin{equation*}
\begin{split}
&\esssup_{t_1<t<t_2}\int_{\Omega}\mathcal{J}((u{-}k)_\pm)\varphi^p\, d\mu
+\int_{t_1}^{t_2}\int_{\Omega} |\nabla (u{-}k)_\pm\varphi|^{p}
\, d\nu \\
&\le C\int_{t_1}^{t_2}\int_{\Omega}(u{-}k)_\pm^p|\nabla \varphi|^p\, d\nu
+C\int_{t_1}^{t_2}\int_{\Omega} \mathcal{J}((u{-}k)_\pm)\varphi^{p-1}\left(\frac{\partial \varphi}{\partial t}\right)_+\, d\nu , \\
\end{split}
\end{equation*}
for every nonnegative $\varphi\in C_0^\infty(\Omega\times(t_1,t_2))$.

\begin{proof}
Let $t_1<\tau_1<\tau_2<t_2$.
We formally substitute the test function $\phi=\pm(u{-}k)_{\pm}\varphi^p$
in the equation and obtain
\begin{equation}\label{substitution}
\begin{split}
0=&\int_{\tau_1}^{\tau_2 }\int_{\Omega} \left( |\nabla u|^{p-2}\nabla u\cdot \nabla \phi
+\frac{\partial (u^{p-1})}{\partial t} \phi \right)\, d\nu  \\
=&\int_{\tau_1}^{\tau_2}\int_{\Omega} |\nabla (u{-}k)|^{p-2}(\pm \nabla (u{-}k)_\pm)
\cdot  \nabla (\pm (u{-}k)_{\pm}\varphi^p )\, d\nu
\\
&\qquad\pm  \int_{\tau_1}^{\tau_2}\int_{\Omega} 
 \frac{\partial (u^{p-1})}{\partial t}(u{-}k)_{\pm}\varphi^p\, d\nu.
\end{split}
\end{equation}
Now the first term on the right-hand side can be estimated pointwise from below as
\begin{align*}
& |\nabla (u{-}k)|^{p-2}(\pm \nabla (u{-}k)_\pm)
\cdot  \nabla (\pm (u{-}k)_{\pm}\varphi^p )
\\ & \qquad \geq |\nabla (u{-}k)_\pm|^p \varphi^p - p |\nabla (u{-}k)_\pm|^{p-1}\varphi^{p-1} (u{-}k)_{\pm} |D\varphi| ,
\end{align*}
and the last term is estimated further by Young's inequality as
\begin{align*}
& - p |\nabla (u{-}k)_\pm|^{p-1}\varphi^{p-1} (u{-}k)_{\pm} |D\varphi| 
\\ & \qquad 
\geq - \frac12 |\nabla (u{-}k)_\pm|^p \varphi^p - C(u{-}k)_{\pm}^p |D\varphi|^p.
\end{align*}%
%
Thus we have
\begin{equation}\label{gradterm_estimate}
\begin{split}
\frac{1}{2}\int_{\tau_1}^{\tau_2}\!&\int_{\Omega} |\nabla (u{-}k)_\pm\varphi|^{p}
\, d\nu\pm\int_{\tau_1}^{\tau_2}\!\int_{\Omega}\frac{\partial (u^{p-1})}{\partial t}(u{-}k)_{\pm}\varphi^p \, d\nu \\
&\le C\int_{\tau_1}^{\tau_2}\int_{\Omega}(u{-}k)_\pm^p|\nabla \varphi|^p\, d\nu.
\end{split}
\end{equation}

Using~\eqref{J_derivative} and integrating by parts,
\begin{align*}
& \pm\int_{\tau_1}^{\tau_2}\int_{\Omega} \frac{\partial (u^{p-1})}{\partial t}(u{-}k)_\pm\varphi^p\, d\nu =\int_{\tau_1}^{\tau_2}\int_{\Omega} \frac{\partial }{\partial t}\mathcal{J}((u{-}k)_\pm)\varphi^p\, d\nu\\
&\qquad = \left[\int_{\Omega}\mathcal{J}((u(x,t){-}k)_\pm)\varphi^p(x,t)\, d\mu\right]_{t=\tau_1}^{\tau_2} 
\\ & \qquad \quad -p\int_{\tau_1}^{\tau_2}\int_{\Omega} \mathcal{J}((u{-}k)_\pm)\varphi^{p-1}\frac{\partial \varphi}{\partial t}\, d\nu.
\end{align*}
So we obtain
\begin{equation}\label{final}
\begin{split}
&\int_{\Omega}\mathcal{J}((u(x,\tau_2){-}k)_\pm)\varphi^p(x,\tau_2)\, d\mu+\int_{\tau_1}^{\tau_2}\int_{\Omega} |\nabla (u{-}k)_\pm\varphi|^{p}
\, d\nu \\
&\le C\int_{\Omega}\mathcal{J}((u(x,\tau_1){-}k)_\pm)\varphi^p(x,\tau_1) \, d\mu+C\int_{\tau_1}^{\tau_2}\int_{\Omega}(u{-}k)_\pm^p|\nabla \varphi|^p\, d\nu \\
&\quad+C\int_{\tau_1}^{\tau_2}\int_{\Omega} \mathcal{J}((u{-}k)_\pm)\varphi^{p-1}\left(\frac{\partial \varphi}{\partial t}\right)_+\, d\nu.
\end{split}
\end{equation}
Now we can drop the second term from the left hand side, let $\tau_1\rightarrow t_1$, choose $\tau_2$ such that
\begin{equation}
\int_\Omega\mathcal{J}((u(x,\tau_2){-}k)_\pm)\varphi^p(x,\tau_2)\, d\mu
\ge\frac{1}{2}\esssup_{t_1<t<t_2}\int_{\Omega}\mathcal{J}((u{-}k)_\pm)\varphi^p\, d\mu
\end{equation}
and estimate the limits of integration on the right hand side of~\eqref{final}. On the other hand, we can also drop the first term on the left hand side of~\eqref{final} and let $\tau_1\rightarrow t_1$ and $\tau_2\rightarrow t_2$. Summing the estimates for both terms gives the claim.
\end{proof}
\end{lemma}

%

Let us denote
\[
\psi_\pm(u):=\Psi(H_k^\pm, (u{-}k)_\pm, c)=\l(\ln\l(\frac{H_k^\pm}{c+H_k^\pm-(u{-}k)_\pm}\r)\r)_+.
\]
The following logarithmic lemma is used in forwarding information in time.

\begin{lemma}\label{Harnack_logarithm}
Let $u\ge 0$ be a weak solution of equation~\eqref{equation}. Then there exists a constant $C=C(p)>0$ such that 
\begin{align*}
\esssup_{t_1<t<t_2}\int_{\Omega} u^{p-2}\psi_-^2(u) \varphi^p \, d\mu
\le & \int_{\Omega} k^{p-2}\psi_-^2(u)(x,t_1)\varphi^p(x) \, d\mu \\
& +C \int_{t_1}^{t_2}\int_\Omega \psi_-|(\psi_-)'|^{2-p}|\nabla  \varphi|^p \, d\nu
\end{align*}
and
\begin{align*}
\esssup_{t_1<t<t_2}\int_{\Omega} k^{p-2}\psi_+^2(u)\varphi^p \, d\mu
\le &\int_{\Omega} u^{p-2}\psi_+^2(u)(x,t_1)\varphi^p(x) \, d\mu 
\\ &+C \int_{t_1}^{t_2}\int_\Omega \psi_+|(\psi_+)'|^{2-p}|\nabla  \varphi|^p \, d\nu.
\end{align*}
Above, $\varphi \in C_0^\infty(\Omega)$ is any nonnegative time-independent test function.

\begin{proof}
Choose
\[
\phi_\pm(u)=\frac{\partial}{\partial u}(\psi_\pm^2(u))\varphi^p
\]
in the definition of weak solution and observe that
\begin{equation}\label{psi''}
(\psi_\pm^2)''=(1+\psi_\pm)(\psi_\pm')^2.
\end{equation}

The parabolic term will take the form
\begin{align*}
&\int_{t_1}^{t_2}\int_\Omega \frac{\partial}{\partial t} u^{p-1} \phi_\pm(u)\, d\nu= \int_{t_1}^{t_2}\int_\Omega \frac{\partial}{\partial t} \int_{k^{p-1}}^{u^{p-1}} \phi_\pm(s^{1/(p-1)}) \, ds\, d\nu \\
&=\l[\int_\Omega\int_{k^{p-1}}^{u^{p-1}} \phi_\pm(s^{1/(p-1)}) \, ds\, d\mu\r]_{t_1}^{t_2} \\
&=\l[(p-1) \int_\Omega\int_{k}^{u} \phi_\pm(r)r^{p-2} \, dr\, d\mu\r]_{t_1}^{t_2}.
\end{align*}
Now an integration by parts gives
\begin{align*}
\int_{k}^{u} \phi_\pm(r)r^{p-2} \, dr
= &\int_{k}^{u} (\psi_\pm^2(r))'r^{p-2} \, dr\varphi^p\\
=&\varphi^p\l[\psi_\pm^2(r)r^{p-2}\r]_{k}^u-(p-2)\int_{k}^{u} \psi_\pm^2(r)r^{p-3} \, dr\varphi^p \\
=&\psi_\pm^2(u)u^{p-2}\varphi^p-(p-2)\int_{k}^{u} \psi_\pm^2(r)r^{p-3} \, dr\varphi^p.
\end{align*}
In the plus case,
we have 
\begin{align*}
\int_{k}^{u} & \phi_+(r)r^{p-2} \, dr
\\ & \ge \psi_+^2(u)u^{p-2}\varphi^p-\psi_+^2(u)(u^{p-2}-k^{p-2})\varphi^p \\
&= (p-1)\psi_+^2(u)k^{p-2}\varphi^p
\end{align*}
and trivially
\[
\int_{k}^{u} \phi_+(r)r^{p-2} \, dr \le \psi_+^2(u)u^{p-2}\varphi^p,
\]
since $p \geq 2$.
Similar estimates are true also for the minus case.
%

On the other hand, by using~\eqref{psi''} together with Young's inequality, we obtain
\begin{align*}
&  |\nabla u|^{p-2}\nabla u \cdot \nabla \phi_\pm
  =  |\nabla u|^{p-2}\nabla u \cdot \nabla((\psi_\pm^2(u))'\varphi^p)
\\ & \qquad \geq |\nabla u|^{p}(1+\psi_\pm)(\psi_\pm')^2\varphi^p
-2p|\nabla u|^{p-1}  \psi_\pm\psi_\pm'\varphi^{p-1}|\nabla\varphi|
\\ & \qquad \geq  \frac12 |\nabla u|^{p}(1+\psi_\pm)(\psi_\pm')^2\varphi^p
- C\psi_\pm(\psi_\pm')^{2-p}|\nabla\varphi|^p,
\end{align*}
almost everywhere, from which the claim follows.
\end{proof}
\end{lemma}

We will need the following notations in the next lemma, which is the most crucial part of the argument. Let
\[
r_n=\frac{r}{2}+\frac{r}{2^{n+1}},
\qquad Q_{n}^\pm=B_{n}\times T_{n}^\pm=B(r_n)\times (t^*-\gamma^\pm r_n^p,t^*)
\]
and
\[
A_{n}^\pm=\left\{(x,t)\in Q_{n}^\pm:\pm u(x,t)>\pm k_n^\pm\right\},
\]
for $n=0,1,2,\dots$. 

Recall the definitions $\mu_i^+-\mu_i^-=\omega_i=\sigma^i\omega_0$, where $\mu_i^+\ge \esssup_{Q^i}{u}$ and $\mu_i^-\le \essinf_{Q^i}{u}$, for $i\ge 1$. Observe, however, that we have to choose 
$$\mu_0^+\le\esssup_K{u} \qquad \textrm{and} \qquad \mu_0^-\ge \essinf_K{u},$$ 
where the infimum and the supremum are taken over $K$ instead of $Q^0$. This is because we need the argument to be independent of the initial cylinder $Q^0$. Now we are ready to prove the fundamental lemma everything depends upon.

\begin{lemma}\label{main_lemma}
Let $0<\eps_\pm \leq 1$, $(k_n^+)_n$ be an increasing sequence and $(k_n^-)_n$ a decreasing sequence, both of nonnegative real numbers. 
Suppose $u\ge0$ is a weak solution of equation~\eqref{equation}, 
\[(u-k_n^\pm)_\pm\le \eps_\pm\omega_i \qquad\text{and}\qquad
|k_{n+1}^\pm-k_n^\pm|\ge\frac{\eps_\pm\omega_i}{2^{n+2}}.
\]
In addition, assume further 
\begin{equation}\label{odd_assumption}
u\ge \frac{1}{C_0}k_n^-  
\end{equation}
and
\begin{equation}\label{nova_assumption}
\mu_i^+\le 2k_n^+, \quad n=1,2,\dots
\end{equation}
for the minus and plus cases, respectively. Then there exist constants \\ $C_-=C(D_0,P_0,C_0,p)>0$ and $C_+=C(D_0.P_0,p)>0$ such that
\begin{equation}\label{minus_iteration}
\frac{\nu(A_{n+1}^\pm)}{\nu(Q_{n+1}^\pm)}\le C_\pm^{n+1} \Gamma_\pm \left(  \frac{\nu(A_n^\pm)}{\nu(Q_n^\pm)}\right)^{2-p/\kappa}
\end{equation}
for every $n=0,1,2,\dots$. Here $\kappa$ is the Sobolev exponent as in \eqref{kappa} and
\[
\Gamma_\pm = \frac{1}{\gamma^\pm}  \left(\frac{k_n^\pm}{\eps_\pm \omega_i}\right)^{p-2} 
\left(\gamma^\pm \left(\frac{\eps_\pm \omega_i}{k_n^\pm}\right)^{p-2} +  1 \right)^{2-p/\kappa}.
\]
\begin{proof}
Choose the cutoff functions $\varphi_n^\pm \in C_0^\infty(Q_n^\pm)$ so that $0\le\varphi_n^\pm\le 1$,
$\varphi_n^\pm=1$ in $Q_{n+1}^\pm$ and
\begin{align}\label{gradient_estimates}
|\nabla \varphi_n^\pm|\le\frac{C2^{n+1}}{r}
\quad\text{and}\quad \left|\frac{\partial \varphi_n^\pm}{\partial t}\right|\le \frac{C2^{p(n+1)}}{\gamma^\pm r^p}.
\end{align}

Denote in short
\[
 \quad v_n = (u{-}k_n)_\pm, \qquad k_n = k_n^\pm, \qquad \eps = \eps_\pm 
\]
and
\[
Q_n = B_n \times T_n = Q_{n}^\pm, \qquad A_n=A_n^\pm,\qquad \gamma=\gamma^\pm, \qquad \varphi_n = \varphi_n^\pm.
\]

By H\"older's inequality, together with the Sobolev inequality \eqref{Sobolevzero}, we obtain
\begin{equation}\label{HoldSobo temp}
 \begin{split}
\dashint_{Q_{n+1}} & v_n^{2(1-p/\kappa)+ p} \, d\nu 
\\ \leq &  
\frac{\nu(Q_{n})}{\nu(Q_{n+1})}
\dashint_{Q_{n}}  v_n^{2(1-p/\kappa)+ p} \varphi_n^{p(1-p/\kappa)+ p} \, d\nu \\
\leq & C \dashint_{T_{n}}\left(\dashint_{B_{n}} v_n^2\varphi_n^p\, d\mu\right)^{1-p/\kappa} 
\left(\dashint_{B_{n}}(v_n \varphi_n)^{\kappa}\, d\mu\right)^{p/\kappa}\, dt \\
\leq & 
C r^p \left(\esssup_{T_n} \dashint_{B_n} v_n^2\varphi_n^p \, d\mu\right)^{1-p/\kappa}
\dashint_{Q_n} |\nabla (v_n\varphi_n)|^p\,d\nu.
\end{split}
\end{equation}
Here, we applied the doubling property of the measure $\nu$ giving 
\[
\frac{\nu(Q_{n})}{\nu(Q_{n+1})} \leq C.
\] 

We continue by studying the term involving the essential supremum.
By the assumption
\[
u\ge \frac{1}{C_0}k_n^-
\]
and~\eqref{minus_lower_estimate2}, we obtain
\begin{align*}
 (u-k_n^-)_-^2 &\leq \frac2{p-1} u^{2-p} \mathcal{J}((u-k_n^-)_-) \leq C (k_n^-)^{2-p}\mathcal{J}((u-k_n^-)_-).
\end{align*}
On the other hand,
the lower bound~\eqref{plus_lower_estimate} gives immediately
\begin{align*}
(u-k_n^+)_+^{2} \leq C (k_n^+)^{2-p} \mathcal{J}((u-k_n^+)_+).
\end{align*}
Using these estimates together with the energy estimate, Lemma~\ref{energy}, yields
\begin{align*}
\esssup_{T_n} & \dashint_{B_n} v_n^2\varphi_n^p \, d\mu
 \leq 
C (k_n)^{2-p} 
\esssup_{T_n} \dashint_{B_n} \mathcal{J}(v_n)\varphi_n^p\, d\mu
\\ \leq & 
C (k_n)^{2-p}  \gamma r_n^p \dashint_{Q_n} \left(  v_n^p |\nabla \varphi_n|^p 
+  \mathcal{J}(v_n) \varphi_n^{p-1}\left(\frac{\partial \varphi_n}{\partial t}\right)_+  \right) \,d\nu. 
\end{align*}
Furthermore, the estimates~\eqref{plus_upper_estimate} and~\eqref{minus_upper_estimate} imply
\[
\mathcal{J}((u{-}k_n)_+) \leq C (k_n^+)^{p-2}(u{-}k_n)_+^2 ,\qquad  \mathcal{J}((u{-}k_n)_-) \leq C (k_n^-)^{p-2}(u{-}k_n)_-^2.
\]
For the plus case, we used \eqref{nova_assumption}. Next, using~\eqref{gradient_estimates}, we arrive at
\begin{equation}\label{eq:main lemma estimate}
\begin{split}
& r_n^p \dashint_{Q_n} \left(  v_n^p |\nabla \varphi_n|^p
+  \mathcal{J}(v_n) \varphi_n^{p-1}\left(\frac{\partial \varphi_n}{\partial t}\right)_+  \right) \,d\nu  \\
& \qquad  \leq
C  2^{np}\dashint_{Q_n} \left(   v_n^p + \frac{(k_n)^{p-2}}{\gamma} v_n^2  \right)\, d\nu \\
& \qquad \leq C 2^{np} (\eps \omega_i)^p
\left(1 + \frac1{\gamma} \left(\frac{\eps \omega_i}{k_n}\right)^{2-p}  \right) \frac{\nu(A_n)}{\nu(Q_n)},
\end{split}
\end{equation}
where the last inequality follows from the fact that $(u{-}k_n)_\pm \leq \eps_\pm \omega_i$.
Thus, we conclude
\begin{equation}\label{eq:main lemma 1}
\esssup_{T_n}  \dashint_{B_n} v_n^2\varphi_n^p \, d\mu \leq 
C 2^{np} (\eps \omega_i)^2
\left(\gamma \left(\frac{\eps \omega_i}{k_n}\right)^{p-2} +  1 \right) \frac{\nu(A_n)}{\nu(Q_n)}.
\end{equation}

Furthermore, since
\[
\dashint_{Q_n} |\nabla (v_n\varphi_n)|^p\,d\nu \leq
C\dashint_{Q_n} |\nabla v_n|^p \varphi_n^p\,d\nu + C\dashint_{Q_n} v_n^p |\nabla\varphi_n|^p\,d\nu, 
\]
applying again the energy estimate and~\eqref{eq:main lemma estimate} leads to
\begin{equation}\label{eq:main lemma 2}
\begin{split}
\dashint_{Q_n} |\nabla (v_n\varphi_n)|^p\,d\nu
\leq &
C 2^{np}  \left(\eps \omega_i\right)^p
\left(1 + \frac1{\gamma} \left(\frac{\eps \omega_i}{k_n}\right)^{2-p}  \right) \frac{\nu(A_n)}{\nu(Q_n)}.
\end{split}
\end{equation}

To finish the proof, note first that
\begin{align*}
(u-k_{n}^\pm)_\pm \chi_{\{(u-k_n^\pm)_\pm>0\}}  \geq &
(u-k_{n}^\pm)_\pm \chi_{\{(u-k_{n+1}^\pm)_\pm>0\}} 
\\ \geq & 
|k_{n+1}^\pm-k_n^\pm| 
\\ \geq  & 2^{-(n+2)}\eps_\pm \omega_i.  
\end{align*}
It then follows that
\begin{equation}\label{eq:main lemma 3}
\dashint_{Q_{n+1}}  v_n^{2(1-p/\kappa)+ p} \, d\nu  \geq  \left(2^{-(n+2)} \eps \omega_i\right)^{2(1-p/\kappa)+ p} 
\frac{\nu(A_{n+1})}{\nu(Q_{n+1})}.
\end{equation}
Inserting estimates~\eqref{eq:main lemma 1}, ~\eqref{eq:main lemma 2}, 
and ~\eqref{eq:main lemma 3} into~\eqref{HoldSobo temp} concludes the proof.
\end{proof}
\end{lemma}

\begin{remark}\label{main_lemma_remark}
 If we have the extra knowledge that $(u-k_n^-)_-=0$, almost everywhere in $B(r)$ at a given time level,
 we can choose the test functions to be independent of time and the cylinder $Q_n^-$ so that
 the length in the time direction stays constant, and the bottom of the cylinder stays
 at the given time level. In this case, by choosing a time independent test function,
 the right hand side of the energy estimate simplifies so that we can get rid of the term $+1$ 
 in the formulation and
\[
\Gamma_\pm = \left(\gamma^\pm \left(\frac{\eps_\pm \omega_i}{k_n^\pm}\right)^{p-2} \right)^{1-p/\kappa}
\]
in~\eqref{minus_iteration}. This will get us the required extra room in the end of the first alternative of Case II.

Furthermore, in the previous lemma, we chose the radii of the cylinder as
\[
 r_n=\frac{r}{2}+\frac{r}{2^{n+1}}.
\]
However, the factor $2$ in the denominator can naturally be replaced by any greater number.

We start the proof by considering Case I. Here, we use the previous lemma only in the plus case. Consequently, the first case does not depend on the constant $C_0$.
\end{remark}

We recall a lemma on the fast geometric convergence of sequences from \cite{DiBe93}.

\begin{lemma}\label{geometric_convergence}
Let $(Y_n)_n$ be a sequence of positive numbers, satisfying
\begin{equation}\label{rec-ineq}
Y_{n+1}\le Cb^nY_n^{1+\alpha}
\end{equation}
where $C,b>1$ and $\alpha>0$. Then $(Y_n)_n$ converges to zero as $n\rightarrow\infty$ provided
\begin{equation}\label{strt-engine}
Y_0\le C^{-1/\alpha}b^{1-\alpha^2}.
\end{equation}
\end{lemma}

On several occasions in the sequel, we use this lemma, together with the fundamental estimate Lemma~\ref{main_lemma}, to conclude that a ratio of the form
\[
Y_n =  \frac{\nu(A_n^\pm)}{\nu(Q_n^\pm)}
\]
converges to zero and consequently that $\nu(A_n^\pm)\rightarrow 0$ as $n\rightarrow \infty$. This will ultimately lead to a reduction of the oscillation which is our final goal. 

Once a recursive inequality of type \eqref{minus_iteration} has been established, the convergence to zero of $\nu(A_n^\pm)$ follows from the condition\[
 \frac{\nu(A_0^\pm)}{\nu(Q_0^\pm)} \le \alpha_0^\pm,
\]
with
\begin{equation}\label{alpha}
\alpha_0^\pm=\Gamma_\pm^{-1/(1-p/\kappa)}   C_\pm^{-1/(1-p/\kappa)+1-(1-p/\kappa)^2},
\end{equation}
where the constants $C_\pm$ and $\Gamma_\pm$ are the constants from the previous lemma. Note that an explicit value of $\alpha_0^\pm$ only follows after fixing $C_\pm$ and $\Gamma_\pm$.

\section{The Case I}

Now we assume that~\eqref{bounds_inf} holds. Our aim is to show that the measures of certain distribution sets
tend to zero and that the local H\"older continuity follows from this.

We start by studying the subcylinder $Q(r^p,r)\subset Q(\eta r^p,r)$. Let $\gamma^\pm=1$, $\eps_\pm=2^{-1}$ and
\[
 k_n^+=\mu_0^+ - \frac{\omega_0}{4}-\frac{\omega_0}{2^{n+2}}.
\]
Observe that, after fixing these quantities, the constant $\alpha_0^+$ can be fixed as well.

We will study two different alternatives which are considered in the following two lemmata, respectively.

\begin{lemma}
Let $\lambda_2>1$ be sufficiently large and let $u\ge 0$ be a weak solution of equation~\eqref{equation}. Furthermore, assume
\begin{equation}\label{measure_assumption2}
\nu\l(\{(x,t)\in B(r)\times (-r^p,-\frac{r^p}{\lambda_2}):u(x,t)\ge \mu_0^-+\frac{\omega_0}{2}\}\r)=0.
\end{equation}
Then there exists a constant $\sigma\in(0,1)$ such that
\[
\essosc_{Q\l(\l(\frac{ r}{2}\r)^p,\frac{r}{2}\r)}{u}\le \sigma \omega_0.
\]
\begin{proof}
By the choices preceding the statement of this lemma, we have
\[
 (u-k_n^+)_+\le \eps_+\omega_0. 
\]
The assumption~\eqref{bounds_inf} implies
\[
\mu_0^+=\mu_0^-+\omega_0\le 2\omega_0.
\]
Thus
\[
1\le\frac{k_n^+}{\eps_+\omega_0}\le 4.
\]
Plug these in Lemma~\ref{main_lemma} to deduce
\[
\frac{\nu(A_{n+1}^+)}{\nu(Q_{n+1}^+)}\le C^{n+1}\left(\frac{\nu(A_n^+)}{\nu(Q_n^+)}\right)^{2-p/\kappa}.
\]

On the other hand, by~\eqref{measure_assumption2} we have the trivial estimate
\[
 \frac{\nu(A_0^+)}{\nu(Q_0^+)}\le \frac{1}{\lambda_2}\le\alpha_0^+,
\]
choosing $\lambda_2>1$ sufficiently large.
By Lemma~\ref{geometric_convergence} we conclude that
\[
 \frac{\nu(A_{n}^+)}{\nu(Q_{n}^+)}\rightarrow 0
\]
as $n\rightarrow \infty$. This implies
\[
\esssup_{Q\left(\left(\frac{r}{2}\right)^p,\frac{r}{2}\right)}{u}
\le \mu_0^+-\frac{\omega_0}{4}.
\]
So, if this alternative occurs, we choose
\[
\mu_1^+=\mu_0^+-\frac{\omega_0}{4}
\]
and
\[
\mu_1^-=\mu_0^-.
\]
These choices yield
\[
\essosc_{Q\left(\left(\frac{r}{2}\right)^p,\frac{r}{2}\right)}{u}\le \left(1-\frac{1}{4}\right)\omega_0
\]
as required, with
\[
\sigma=\frac{3}{4}.
\]
\end{proof}

\end{lemma}

For the second possibility, we have the following lemma.

\begin{lemma}\label{forwarding}
Let $u\ge 0$ be a weak solution of equation~\eqref{equation} and suppose
\begin{equation}\label{measure_assumption}
\nu\l(\{(x,t)\in B(r)\times (-r^p,-\frac{r^p}{\lambda_2}):u(x,t)\ge \mu_0^-+\frac{\omega_0}{2}\}\r)>0.
\end{equation}
Then there exists a constant $\sigma=\sigma(H_0)\in(0,1)$ such that
\[
\essosc_{Q\left(\left(\frac{r}{2\lambda_2}\right)^p, \frac{r}{2\lambda_2}\right)}{u}\le \sigma\omega_0.
\]

\begin{proof}
By assumption~\eqref{measure_assumption}, we have
\[
 \esssup_{B(r)\times(-r^p,-\frac{r^p}{\lambda_2})}{u}\ge \mu_0^-+\frac{\omega_0}{2}.
\]
Now we can use Harnack's inequality (Theorem~\ref{Harnack}), together with the Case I assumption~\eqref{bounds_inf}, to deduce
\begin{align*}
\essinf_{Q\left(\left(\frac{r}{2\lambda_2}\right)^p,\frac{r}{2\lambda_2}\right)}{u}&\ge\frac{1}{H_0}\esssup_{B(r)\times(-r^p,-\frac{r^p}{\lambda_2})}{u} \\
&\ge  \frac{\mu_0^-}{H_0}+\frac{\omega_0}{2H_0} \\
&\ge \mu_0^-+ \frac{\mu_0^-}{H_0}-\frac{\omega_0}{2H_0+1}+\frac{\omega_0}{2H_0} \\
&\ge \mu_0^-+\frac{\omega_0}{2H_0(2H_0+1)}.
\end{align*}
Observe that the constant $H_0$ depends on $\lambda_2$, but this does not matter since $\lambda_2$ depends only on the data. 

Now, if we end up in this alternative, we choose
\[
\mu_1^-=\mu_0^-+\frac{\omega_0}{2H_0(2H_0+1)}
\]
and
\[
\mu_1^+=\mu_0^+.
\]
We also obtain
\[
\essosc_{Q\left(\left(\frac{r}{2\lambda_2}\right)^p,\frac{r}{2\lambda_2}\right)}{u}
\le \omega_0-\frac{\omega_0}{2H_0(2H_0+1)}
=\sigma\omega_0,
\]
with
\[
\sigma=1-\frac{1}{2H_0(2H_0+1)},
\]
as required.
\end{proof}
\end{lemma}

\section{The case II}

In Case II the equation looks like the evolution $p$-Laplace equation. In this case, we need to use the scaling factor $\eta$ in the time geometry of our cylinders. The difficulty is now that we cannot use the Harnack principle anymore, as the lower bound it gives might be trivial. Indeed, the infimum can be larger than the lower bound Harnack's inequality gives. On the other hand, we have the following kind of elliptic Harnack's inequality.

Suppose that $j_0$ is the first index for which assumption~\eqref{bounds_inf} does not hold.
Then we have
\begin{equation}\label{Harnack_assumption}
 \omega_{j_0}\le\mu_{j_0}^+\le (2H_0+2)\mu_{j_0}^-.
\end{equation}
Clearly, this Harnack's inequality is valid also for every subset of the initial cylinder $Q_{j_0}=Q(\eta r^p,r)$ and, consequently, for every $j\ge j_0$.

Recall, that $\omega_{j_0}=\sigma\omega_{j_0-1}$ and
\begin{equation}
\begin{split}
\frac{\omega_{j_0}}{(2H_0+2)} &\le\mu_{j_0}^- \le\mu_{j_0-1}^-+(1-\sigma)\omega_{j_0-1} \\
&\le(2-\sigma)\omega_{j_0-1} \le\frac{2-\sigma}{\sigma}\omega_{j_0}.
\end{split}
\end{equation}
Thus, we obtain
\[
\frac{\sigma}{(2H_0+2)(2-\sigma)}\le\frac{\sigma}{2-\sigma}\frac{\mu_{j_0}^-}{\omega_{j_0}}\le 1
\]
and, consequently,
\[
Q\l(\frac{C_1}{(2H_0+2)^{p-2}}\eta r^p,r\r) \subset Q\left(C_1 \left(\frac{\mu_{j_0}^-}{\omega_{j_0}}\right)^{p-2}\eta r^p,r\right) \subset Q(\eta r^p,r),
\]
where $C_1=\sigma^{p-2}/(2-\sigma)^{p-2}$. We will consider the cylinder
\[
Q:=Q\left(C_1 \left(\frac{\mu_{j_0}^-}{\omega_{j_0}}\right)^{p-2}\eta r^p,r\right).
\]
By the above calculation, we have shrunk the cylinder by a factor which is controllable by the data.

In the sequel, we will denote
\begin{equation}\label{theta}
\theta=C_1\left(\frac{\mu_{j_0}^-}{\omega_{j_0}}\right)^{p-2}.
\end{equation}

Recall the definitions
\[
 Q_n^\pm=Q_{t^*}(\gamma^\pm r_n^p,r_n)=B_n\times T_n=B(r_n)\times (t^*-\gamma^\pm r_n^p,t^*)
\]
and
\[
A_n^\pm=\{(x,t)\in Q_n^\pm:\pm u>\pm k_n^\pm\}.
\]

Now the proof will follow the classical argument of DiBenedetto, see~\cite{DiBe93} and~\cite{Urba08}, and is again divided into two alternatives. In the first one, we assume that there is a suitable cylinder for which the set where $u$ is close to its infimum is very small. In the second alternative, we assume that this does not hold true.

\subsection{The First Alternative}

We first  suppose that there exists a constant $\alpha_0 \in (0,1)$ (to be determined in the course of the next lemma, depending only on the data) such that 

$$\nu \left( \{ (x,t)     \in Q_0^- : u < \mu_{j_0}^- + \frac{\omega_{j_0}}{2} \} \right) \leq     \alpha_0 \nu ( Q_0^-),$$
for a cylinder
\[
Q_0^-=Q_{t^*}(\theta r^p,r)\subset Q(2^{\lambda_1(p-2)+1}\theta r^p,r).
\]

Our aim is to use Lemma~\ref{main_lemma} to conclude for the reduction of the oscillation. 

\begin{lemma}\label{time_forwarding}
For every $s>3$,
\begin{align*}
&\nu\left(\left\{(x,t)\in Q\left(\theta\l(\frac{r}{4}\r)^p,\frac{r}{4}\right):u(x,t)<\mu_{j_0}^-+\frac{\omega_{j_0}}{2^{s}}
\right\}\right)
\\&\qquad\qquad
\le C2^{\lambda_1(p-2)}\frac{s-2}{(s-3)^2}\nu\left(Q\left(\theta\l(\frac{r}{4}\r)^p,\frac{r}{4}\right)\right),
\end{align*}
where 
$\theta$ is as in~\eqref{theta}.

\begin{proof}
We start by using Lemma~\ref{main_lemma}, with the choices
\[
r_n=\frac{r}{2}+\frac{r}{2^{n+1}},\qquad k_n^-=\mu_{j_0}^-+\frac{\omega_{j_0}}{4}+\frac{\omega_{j_0}}{2^{n+2}},
\]
$\eps_-=1/2$ and $\gamma^-=\theta$. We also need the assumption~\eqref{Harnack_assumption} to deduce that
\begin{equation}\label{C_0}
u\ge \frac{1}{C_0}(\mu_{j_0}^-+2(2H_0+2)\mu_{j_0}^-)\ge\frac1{C_0}(\mu_{j_0}^-+2\omega_{j_0})\ge \frac1{C_0}k_n^-,
\end{equation}
with $C_0=3(2H_0+2)$.  This knowledge is needed in Lemma~\ref{main_lemma}. Now, after fixing $\eps_-$, $\gamma^-$, $k_n^-$ and $C_0$ we can fix $\alpha_0^-$, see~\eqref{alpha}. 

We also obtain, using \eqref{Harnack_assumption} and \eqref{theta}, the bounds
\[
\frac{2^{p-2}}{C_1}\le \frac{1}{\gamma^-}\l(\frac{k_n^-}{\eps_- \omega_{j_0}}\r)^{p-2}\le\frac{2^{p-2}}{C_1}(2H_0+2)^{p-2}
\]
and thus we can conclude
\[
\frac{\nu(A_{n+1}^-)}{\nu(Q_{n+1}^-)}\le C^{n+1}\left(\frac{\nu(A_n^-)}{\nu(Q_n^-)}\right)^{2-p/\kappa}.
\]

By the assumption of this alternative, together with the lemma of fast geometric convergence (Lemma~\ref{geometric_convergence}), we have $u>k$ almost everywhere in $Q_{t^*}(\theta (r/2)^p,r/2)$. Thus
\[
(u{-}k)_-=0
\]
and consequently 
\[
\psi(u):=\l(\ln\l(\frac{H_k^-}{c+H_k^--(u{-}k)_-}\r)\r)_+=0
\]
almost everywhere in $Q_{t^*}(\theta (r/2)^p,r/2)$. Let
\begin{equation}\label{t_0}
t_0\le -\theta\l(\frac{r}{4}\r)^p
\end{equation}
be a time level such that this is true for almost every $x\in B(r/2)$.

Now our goal is to apply Lemma~\ref{Harnack_logarithm} with
\[
k=\mu_{j_0}^-+\frac{\omega_{j_0}}{4}, \quad c=\frac{\omega_{j_0}}{2^{s}}
\]
and
\[
H_k^-=\esssup_{Q}{(u-k)_-},
\]
where $Q=Q(\eta\theta r^p,r)$. Choose $\varphi\in C_0^\infty(B(r/2))$ independent of time such that
$0\le\varphi\le 1$, $\varphi=1$ in $B(r/4)$
and
\[
|\nabla \varphi|\le\frac{C}{r}.
\]
In the set $\{u< \mu_{j_0}^-+\frac{\omega_{j_0}}{2^{s}}\}$, we have
\[
\psi^2\ge (s-3)^2 \ln^22,
\]
and, on the other hand,
\[
\psi\le (s-2)\ln 2\quad\text{and}\quad|\psi'|^{2-p}\le\l(\frac{\omega_{j_0}}{2}\r)^{p-2}.
\]
The use of these estimates in Lemma~\ref{Harnack_logarithm} gives
\begin{align*}
&(\mu_{j_0}^-)^{p-2}(s-3)^2 \ln^22\cdot\mu(\{x\in B(r/4):u(x,t)<\mu_{j_0}^-+\frac{\omega_{j_0}}{2^{s}}\}) \\
&\le\esssup_{t_0<t<0}\int_{B(r/2)} u^{p-2}\psi^2(u)(x,t)\varphi^p(x) \, d\mu \\
&\le \int_{B(r/2)} k^{p-2}\psi^2(u)(x,t_0)\varphi^p(x) \, d\mu \\
&+C \int_{t_0}^{0}\int_{B(r/2)} \psi|\psi'|^{2-p}|\nabla  \varphi|^p \, d\mu\, dt \\
&\le C(s-2)\ln 2\l(\frac{2^{\lambda_1}\omega_{j_0}}{2}\r)^{p-2}\theta\mu(B(r/4)) \\
&\le C(s-2)\ln 2\l(\frac{2^{\lambda_1}\mu_{j_0}^-}{2}\r)^{p-2}\mu(B(r/4)),
\end{align*}
for almost every $t\in (t_0,0)$. Observe that, in the third inequality, we plugged in $\eta=2^{\lambda_1(p-2)+1}$. The claim follows by integrating this estimate over $(-\theta (r/4)^p,0)$.

\end{proof}
\end{lemma}

We conclude this alternative with the following two lemmata.

\begin{lemma}
Let $u\ge 0$ be a weak solution of equation~\eqref{equation} and assume~\eqref{Harnack_assumption} holds. Then
\[
u\ge \mu_{j_0}^-+\frac{\omega_{j_0}}{2^{s+1}}\qquad\text{a.e. in}\quad Q\left(\theta\left(\frac{r}{8}\right)^p,\frac{r}{8}\right),
\]
where $s$ depends only upon the data, and $\theta$ is as in~\eqref{theta}.

\begin{proof}
Let
\[
r_n=\frac{r}{8}+\frac{r}{2^{n+3}},
\]
\[
 Q_n^-=B_n\times T=B(r_n)\times (t_0,0),
\]
where $t_0$ is, as in the previous lemma, such that
\[
(u{-}k)_-(x,t_0)=0,\quad t_0\le -\theta\l(\frac{r}{4}\r)^p,
\]
for a.e. $x\in B(r/2)$. Moreover, define
\[
k_n^-=\mu_{j_0}^-+\frac{\omega_{j_0}}{2^{s+1}}+\frac{\omega_{j_0}}{2^{s+n+1}}.
\]

In this case, we obtain
\[
 (u-k_n^-)_-\le \eps_-\omega_{j_0}, \quad\text{where}\quad\eps_-=\frac{1}{2^s}.
\]
Observe also that $\gamma^-=-t_0/r^p\le \eta\theta=2^{\lambda_1(p-2)+1}\theta$.

We will substitute these in Lemma~\ref{main_lemma} and, taking into account Remark~\ref{main_lemma_remark} and estimate~\eqref{C_0}, we conclude as before that
\[
\frac{\nu(A_{n+1}^-)}{\nu(Q_{n+1}^-)}\le C^{n+1}\l(\frac{2^{\lambda_1(p-2)}}{2^{s(p-2)}}\r)^{1-p/\kappa}\left(\frac{\nu(A_n^-)}{\nu(Q_n^-)}\right)^{2-p/\kappa}.
\]
Now choose $s>\lambda_1$. Then, by Lemma~\ref{geometric_convergence}, we have $\nu(A_n^-)/\nu(Q_n^-)\rightarrow 0$ as $n\rightarrow \infty$, provided $\nu(A_0^-)/\nu(Q_0^-)$ is small enough.
On the other hand, by choosing $s$ large enough, Lemma~\ref{time_forwarding} guarantees that $\nu(A_0^-)/\nu(Q_0^-)$
can be chosen to be as small as we please.

This gives
\[
u\ge \mu_{j_0}^-+\frac{\omega_{j_0}}{2^{s+1}}\qquad\text{a.e. in}\quad Q\left(|t_0|,\frac{r}{8}\right)
\]
and hence the lemma is proved.

\end{proof}
\end{lemma}

\begin{lemma}
There exists $0<\sigma<1$, depending only upon the data, such that
\[
\essosc_{Q\left(\theta\left(\frac{r}{8}\right)^p,\frac{r}{8}\right)}{u}\le \sigma \omega_{j_0}.
\]

\begin{proof}
By the previous lemma,
\[
\essinf_{Q\left(\theta\left(\frac{r}{8}\right)^p,\frac{r}{8}\right)}{u}\ge \mu_{j_0}^-+\frac{\omega_{j_0}}{2^{s+1}},
\]
for some $s>1$, which depends only upon the data  and $\lambda_1$. Observe that here we used the knowledge
\[
t_0\le-\theta\l(\frac{r}{8}\r)^p.
\]

If this alternative occurs, we again choose
\[
\mu_{j_0+1}^-:=\mu_{j_0}^-+\frac{\omega_{j_0}}{2^{s+1}}
\]
and
\[
\mu_{j_0+1}^+=\mu_{j_0}^+.
\]
Finally, we get
\[
\essosc_{Q\left(\theta\left(\frac{r}{8}\right)^p,\frac{r}{8}\right)}{u}\le \left(1-\frac{1}{2^{s+1}}\right)\omega_{j_0}
\]
as required, with
\[
 \sigma=1-\frac{1}{2^{s+1}}.
\]

\end{proof}

\end{lemma}

\begin{remark}
 Here the choice of $s$ is possible only after $\lambda_1$ has been determined in the second alternative. Nevertheless, both of them are a priori constants which can be assigned explicit values depending only upon the data.
\end{remark}

%
%
%
%
%
%

\subsection{The Second Alternative}

In the second alternative, the assumption of the first alternative is not true.
In this case, for every cylinder $Q_{t^*}(\theta r^p,r) \subset Q(\eta \theta r^p,r)$,
we have
\begin{equation}\label{assumption2}
\begin{split}
\frac{\nu \left( \{ (x,t) \in Q_{t^*} (\theta r^p,r) : u (x,t) \geq  
\mu_{j_0}^- + \frac{\omega_{j_0}}{2} \} \right)}{\nu (  
Q_{t^*} (\theta r^p,r) )} <  (1-\alpha_0) ,
\end{split}
\end{equation}
where $\alpha_0:=\alpha_0^-$ is the same constant as in the first alternative.
This implies that, for every $t^*\in (-(\eta-1)\theta r^p,0)$, there exists a time level $t_0$ with
\[
t^*- \theta r^p\le t_0 \le t^*-\frac{\theta\alpha_0}{2} r^p
\]
for which
\begin{equation}\label{second_alternative}
\begin{split}
\mu \left(\left\{ x\in B(r) : u(x,t_0)>k_0^-\right\}\right)\le \frac{1-\alpha_0}{1-\displaystyle{\frac{\alpha_0}2}}\mu(B(r)).
\end{split}
\end{equation}
Indeed, otherwise we would have
\begin{align*}
&\nu\left(\left\{ (x,t)\in Q_{t^*}(\theta r^p,r) : u(x,t)>k_0^-\right\}\right) \\
&\ge\int_{t^*- \theta r^p}^{t^*-\frac{\theta\alpha_0}{2} r^p}\mu\left(\left\{ x\in B(r) : u(x,t) >k_0^-\right\}\right) \, dt \\
&>(1-\alpha_0)\nu(Q_{t^*}(\theta r^p,r)),
\end{align*}
which contradicts~\eqref{assumption2}.

This alternative is also based on Lemma~\ref{main_lemma}. We choose $\lambda_1$ in the definition of $k_n^+$ large enough so that
we can force $\nu(A_0^+)$ to be small compared to $\nu(Q_0^+)$.


We start with forwarding the information of~\eqref{second_alternative} in time.

\begin{lemma}\label{logarithmic_bound}
There exists $s^*>0$, depending only upon the data, such that
\[
\mu\left(\left\{x\in B(r): u(x,t)>\mu_{j_0}^+-\frac{\omega_{j_0}}{2^{s^*}}
\right\}\right)
\le\frac{1-\displaystyle{\frac{3\alpha_0}4}}{1-\displaystyle{\frac{\alpha_0}2}}\mu(B(r)).
\]
for almost all $t\in\left(t_0,0\right)$.

\begin{proof}
Let
\[
c=\frac{\omega_{j_0}}{2^{s+n}}, \quad k=\mu_{j_0}^+-\frac{\omega_{j_0}}{2^s}
\]
and
\[
H_k^+=\esssup_{Q}{(u{-}k)_+},
\]
where $s$ and $n$ will be chosen later and $Q:=Q(\theta\eta r^p,r)$. Our aim is again to use Lemma~\ref{Harnack_logarithm} to forward the information in time. We will need some estimates for doing this.

Recall the definition
\[
\psi_+(u)=\Psi(H_k^+,(u{-}k)_+,c)=\ln^+\l(\frac{H_k^+}{c+H_k^+-(u{-}k)_+}\r).
\]
Trivially, we have
\[
\psi_+(u)\le \ln\l(\frac{\frac{\omega_{j_0}}{2^s}}{\frac{\omega_{j_0}}{2^{s+n}}}\r)=n\ln 2
\]
and, on the other hand, in the set
\[
 \{u>l \equiv \mu_{j_0}^+-\frac{\omega_{j_0}}{2^{s+n}}\},
\]
we get
\[
\psi_+(u)\ge \ln\l(\frac{\frac{\omega_{j_0}}{2^s}}{\frac{\omega_{j_0}}{2^{s+n}}+\frac{\omega_{j_0}}{2^{s+n}}}\r)=(n-1) \ln 2
\]

The last estimate we need is 
\[
|(\psi_+)'|^{2-p}\le \l(\frac{1}{c+H_k^+}\r)^{2-p}\le 2^{p-2}\l(\frac{\omega_{j_0}}{2^s}\r)^{p-2}.
\]

Let now $\varphi\in C_0^\infty(B(r))$ be a cutoff function which is independent of time
and has the properties $0\le\varphi\le 1$, $\varphi=1$ in $B((1-\delta )r)$ and
\[
|\nabla \varphi|\le\frac{1}{\delta r},
\]
where $0<\delta<1$ is to be determined later.

Apply Lemma~\ref{Harnack_logarithm} with these choices to conclude
\begin{equation}\notag
\begin{split}
&(n-1)^2\ln^22 \cdot\mu(\{x\in B((1-\delta)r):u(x,t)>l\}) \\
&\le\esssup_{t_0<t<t^*}\int_{B(r)} \psi_+^2(u)(x,t)\varphi^p(x) \, d\mu \\
&\le \int_{B(r)} \l(\frac{u}{k}\r)^{p-2}\psi_+^2(u)(x,t_0)\varphi^p(x) \, d\mu \\
&+C k^{2-p}\int_{t_0}^{t^*}\int_{B(r)}\psi_+|(\psi_+)'|^{2-p}|\nabla  \varphi|^p \, d\mu\, dt \\
&\le n^2 \ln^2(2)\l(\frac{\mu_{j_0}^+}{\mu_{j_0}^+-\frac{\mu_{j_0}^+}{2^{s}}}\r)^{p-2} \frac{1-\alpha_0}{1-\displaystyle{\frac{\alpha_0}2}}\mu(B(r)\times\{t_0\}) \\
&\quad+C\frac{n\ln 2}{\delta^p} \l(\frac{\omega_{j_0}}{\mu_{j_0}^+}\r)^{p-2}\l(\frac{1}{2^s-1}\r)^{p-2} \theta\mu(B(r))
\end{split}
\end{equation}
for almost every $t\in(t_0,t^*)$. Observe, that in the third inequality we used~\eqref{second_alternative}.

Now, by the annular decay property \eqref{annular_decay}, we have
\begin{align*}
\mu(\{x\in &B(r):u(x,t)>l\}) \\
&\le\mu(B(r)\setminus B((1-\delta)r))+\mu(\{x \in B((1-\delta) r):u(x,t)>l\})\\
&\le C\delta^\alpha\mu(B(r))+\mu(\{x \in B((1-\delta) r):u(x,t)>l\}).
\end{align*}
For the first term, we choose $\delta$ small enough so that
\[
C\delta^\alpha<\frac{\alpha_0}{16\left(1-\displaystyle{\frac{\alpha_0}2}\right)}
\]
and for the second term we use the previous estimate. Indeed, by choosing $s$ and $n$ large enough so that
\[
\frac{1-\alpha_0}{1-\displaystyle{\frac{\alpha_0}2}}\frac{n^2}{(n-1)^2}\l(\frac{1}{1-\frac{1}{2^{s}}}\r)^{p-2}
\le 1-\frac{3\alpha_0}{8\left(1-\displaystyle{\frac{\alpha_0}2}\right)}
\]
and
\[
\frac{Cn}{(\ln 2)\delta^p(n-1)^2}\l(\frac{1}{2^s-1}\r)^{p-2}
\le\frac{\alpha_0}{16\left(1-\displaystyle{\frac{\alpha_0}2}\right)},
\]
we get the claim for almost every $t\in(t_0,t^*)$ with $s^*=s+n$.
Recall that
\[
t^*- \theta r^p\le t_0\le t^*-\frac{\theta\alpha_0}{2}r^p.
\]
Finally, since the above holds for every cylinder $Q_{t^*}(\theta r^p,r) \subset Q(\eta\theta r^p,r)$, we can conclude that for almost every
\[
t\ge -\eta \theta r^p+\theta r^p-\frac{\theta\alpha_0 r^p}{2}=\left(-\eta+1-\frac{\alpha_0}{2} \right)\theta r^p,
\]
we have
\[
\mu(\{x\in B(r):u(x,t)>l\})<\frac{1-\displaystyle{\frac{3\alpha_0}4}}{1-\displaystyle{\frac{\alpha_0}2}}\mu(B(r)).
\]
\end{proof}
\end{lemma}

\begin{remark}
Now we can choose $\eta$ so large that the previous
lemma holds for almost every $t\in\big(-\frac{\eta\theta}{2}r^p,0\big)$, i.e.,
\[
-\eta+1-\frac{\alpha_0}{2}\le-\frac{\eta}{2}
\]
and hence
\[
\eta=2^{\lambda_1(p-2)+1} \ge 2-\alpha_0.
\]
But this is always guaranteed for $\lambda_1>1$ and $p\ge2$.
\end{remark}

We are ready to prove the final estimate, which, together with
Lemma~\ref{geometric_convergence}, gives the reduction of the oscillation. Let
\[
E_\varrho (t) = \{ x \in B(r) : u(x,t) > \mu_{j_0}^+ -  
\frac{\omega_{j_0}}{2^\varrho}  \}
\]
and
\[
E_\varrho = \{ (x,t) \in Q\l(\frac{\eta \theta}{2} r^p,r\r) : u(x,t)  
 >\mu_{j_0}^+ - \frac{\omega_{j_0}}{2^\varrho}  \}.
\]
Then we have the following lemma.
\begin{lemma}\label{choosing_lambda}
For every $\alpha_1\in(0,1)$, there exists $\lambda_1>0$ such that
\[
\frac{\nu(E_{\lambda_1})}{\nu(Q(\frac{\eta\theta}{2} r^p,r))}\le \alpha_1.
\]

\begin{proof}
Denote
\[
h=\mu_{j_0}^+-\frac{\omega_{j_0}}{2^{s+1}}
\]
and
\[
k=\mu_{j_0}^+-\frac{\omega_{j_0}}{2^{s}},
\]
where $s>0$ will be chosen large. Let also
\[
v=
\begin{cases}
h-k,\quad &u\ge h, \\
u-k,\quad &k<u<h, \\
0, \quad &u\le k.
\end{cases}
\]
By the previous lemma, we can choose $s$ large enough, namely $s\ge s^*$, so that, for
almost every $t\in\big(-\frac{\eta\theta}{2}r^p,0\big)$, we have
\begin{align*}
\mu(x\in B(r): v(x,t)=0\})&=\mu(\{x\in B(r): u(x,t)\le k\}) \\
&\ge\frac{\alpha_0}{4-2\alpha_0}\mu(B(r)\times\{t\}) \\
&\ge\frac{\alpha_0}{4}\mu(B(r)).
\end{align*}
Thus, for almost every $t\in\big(-\frac{\eta\theta}{2}r^p,0\big)$, we obtain
\[
v_{B(r)}(t)=\dashint_{B(r)\times\{t\}} v \, d\mu \le \l(1-\frac{\alpha_0}{4}\r)(h-k)
\]
and, consequently,
\[
h-k-v_{B(r)}(t) \ge \frac{\alpha_0}{4}(h-k).
\]

Using the $(q,q)$-Poincar\'e inequality for some $q<p$ (see~\eqref{poincare} and the remark after that), yields
\begin{align*}
(h-k)^q\mu(E_{s+1}(t))&\le\l(\frac{4}{\alpha_0}\r)^q\int_{B(r)\times\{t\}}|v-v_{B(r)}(t)|^q\, d\mu \\
&\le
Cr^q\int_{B(r)\times\{t\}} |\nabla v|^q \, d\mu = Cr^q\int_{E_s(t)\setminus E_{s+1}(t)} |\nabla u|^q \, d\mu,
\end{align*}
for almost every $t\in\big(-\frac{\eta\theta}{2}r^p,0\big)$. The constant $(4/\alpha_0)^q$ above was absorbed into the constant $C$. Now we integrate the above inequality over time to get
\[
(h-k)^q\nu(E_{s+1})\le C r^q\int_{E_s\setminus E_{s+1}} |\nabla u|^q \, d\nu.
\]
Next, we introduce a cutoff function $\varphi \in C_0^\infty(Q(\eta\theta r^p,2r))$ such that
$0\le\varphi\le 1$, $\varphi=1$ in $Q\left(\frac{\eta\theta}{2} r^p,r\right)$
and
\[
|\nabla \varphi|\le\frac{C}{r} \quad\text{and}\quad \left|\frac{\partial \varphi}{\partial t}\right|\le \frac{C}{\eta \theta r^p}.
\]

Now H\"older's inequality gives
\begin{align*}
(h-k&)^q\nu(E_{s+1})\le Cr^q\l(\int_{E_s\setminus E_{s+1}} |\nabla u|^p \, d\nu\r)^{q/p}\nu(E_s\setminus E_{s+1})^{1-q/p} \\
&\le Cr^q\l(\int_{Q(\eta\theta r^p,2r)} |\nabla (u{-}k)_+|^p\varphi^p \, d\nu\r)^{q/p}\nu(E_s\setminus E_{s+1})^{1-q/p}. \\
\end{align*}

By choosing $\lambda_1>s\ge s^*$ in the definition of $\eta$ large enough,
the first factor on the right hand side can be estimated by Lemma~\ref{energy} and~\eqref{plus_upper_estimate} as
\begin{equation}\label{nabla_u}
\begin{split}
&\int_{Q(\eta\theta r^p,2r)} |\nabla (u{-}k)_+|^p\varphi^p \, d\nu \\
&\le C\int_{Q(\eta \theta r^p,2r)}(u{-}k)_+^p|\nabla\varphi|^p\, d\nu \\
&\quad+C\int_{Q(\eta \theta r^p,2r)} \mathcal{J}((u{-}k)_+)\varphi^{p-1}\left|\frac{\partial \varphi}{\partial t}\right|\, d\nu \\
&\le \frac{C}{r^p}\l(\l(\frac{\omega_{j_0}}{2^s}\r)^{p-2}+\frac{(\mu_{j_0}^+)^{p-2}}{\eta\theta }\r)\int_{Q(\eta \theta r^p,2r)}(u{-}k)_+^2 \, d\nu \\
&\le \frac{C}{ r^p}\left(\frac{\omega_{j_0}}{2^s}\right)^p\nu\left(Q\left( \frac{\eta \theta}{2} r^p,r\right)\right).
\end{split}
\end{equation}
In the last inequality we used the doubling property of the measure $\nu$.

We obtain
\[
\l(\frac{\omega_{j_0}}{2^{s+1}}\r)^q
\nu(E_{s+1})\le C
\l(\frac{\omega_{j_0}}{2^s}\r)^q
\nu\left(Q\left( \frac{\eta \theta}{2} r^p,r\right)\right)^{q/p}\nu(E_s\setminus E_{s+1})^{1-q/p}.
\]
Finally, summing $s$ over $s^*,\dots,\lambda_1-1$ gives
\[
(\lambda_1-s^*)\nu(E_{\lambda_1})^{p/(p-q)}
\le C\nu\left(Q\left( \frac{\eta \theta}{2} r^p,r\right)\right)^{q/(p-q)}\nu\left(Q\left( \frac{\eta \theta}{2} r^p,r\right)\right)
\]
and hence
\[
\nu(E_{\lambda_1})\le \frac{C}{(\lambda_1-s^*)^{(p-q)/p}}\nu\left(Q\left( \frac{\eta \theta}{2} r^p,r\right)\right).
\]
Choosing $\lambda_1$ large enough finishes the proof.
\end{proof}
\end{lemma}

\begin{lemma}
Suppose that~\eqref{assumption2} holds. Then there exists $0<\sigma<1$, depending only upon the data, such that
\[
\essosc_{Q(\frac{\eta\theta}{2}r^p,r)}{u}\le \sigma \omega_{j_0}.
\]

\begin{proof}
Let
\[
 Q_n^+=B(r_n)\times (-\gamma^+ r_n^p,0),
\]
\vspace{0in}
\[
 r_n=\frac{r}{2}+\frac{r}{2^{n+1}}
\]
and $A_n^+$ as before. Substituting $\gamma^+=2^{\lambda_1(p-2)}\theta$, $\eps_+=1/2^{\lambda_1}$ and
\[
 k_n^+=\mu_{j_0}^+ - \frac{\omega_{j_0}}{2^{\lambda_1+1}}-\frac{\omega_{j_0}}{2^{\lambda_1+n+1}}
\]
in Lemma~\ref{main_lemma}, and using~\eqref{Harnack_assumption} to bound
\[
1\le\frac{1}{\gamma^+}\l(\frac{k_n^+}{\eps_+\omega_{j_0}}\r)^{p-2}\le (2H_0+2)^{p-2},
\]
yields
\[
\frac{\nu(A_{n+1}^+)}{\nu(Q_{n+1}^+)}\le C^{n+1}\left(\frac{\nu(A_n^+)}{\nu(Q_n^+)}\right)^{2-p/\kappa}.
\]
By the previous Lemma, we can choose $\lambda_1$ large enough so that
\[
\frac{\nu(A_{0}^+)}{\nu(Q_{0}^+)}
\]
is as small as we please. Consequently, by Lemma~\ref{geometric_convergence}, we obtain
\[
\esssup_{Q\left(\frac{\eta\theta}{2}\left(\frac{r}{2}\right)^p,\frac{r}{2}\right)}{u}
\le \mu_{j_0}^+-\frac{\omega_{j_0}}{2^{\lambda_1+1}},
\]
for some $\lambda_1>1$, which depends only upon the data. So if this alternative occurs, we choose
\[
\mu_1^+=\mu_{j_0}^+-\frac{\omega_{j_0}}{2^{\lambda_1+1}}
\]
and
\[
\mu_1^-=\mu_{j_0}^-.
\]
We also obtain
\[
\essosc_{Q\left(\frac{\eta\theta}{2}\left(\frac{r}{2}\right)^p,\frac{r}{2}\right)}{u}\le \left(1-\frac{1}{2^{\lambda_1+1}}\right)\omega_{j_0},
\]
as required, with
\[
\sigma=1-\frac{1}{2^{\lambda_1+1}}.
\]
\end{proof}

\end{lemma}

Now the H\"older continuity follows by standard iterative real analysis methods, see chapter III of~\cite{DiBe93}, or~\cite{Urba08}.

\bigskip

\noindent {\bf Acknowledgments.} The authors would like to thank Juha Kinnunen for interesting discussions on the subject of this paper.

\bigskip

\bibliography{citations}
\bibliographystyle{plain}

\bigskip
\noindent Addresses:

\noindent T.K.: Aalto University, Institute of Mathematics, P.O. Box 11100, FI-00076 Aalto, Finland. \\
\noindent
E-mail: {\tt tuomo.kuusi@tkk.fi}\\

\noindent J.S.: Aalto University, Institute of Mathematics, P.O. Box 11100, FI-00076 Aalto, Finland. \\
\noindent
E-mail: {\tt juhana.siljander@tkk.fi}\\

\noindent J.M.U.: CMUC, Department of Mathematics,
University of Coimbra, 3001-454 Coimbra, Portugal. \\
\noindent
E-mail: {\tt jmurb@mat.uc.pt}\\

\end{document}